# Embedding 5-planar graphs in three pages

Xiaxia Guan, Weihua Yang*

Department of Mathematics, Taiyuan University of Technology,

Taiyuan Shanxi-030024, China

**Abstract.** A *book-embedding* of a graph $G$ is an embedding of vertices of $G$ along the spine of a book, and edges of $G$ on the pages so that no two edges on the same page intersect. the minimum number of pages in which a graph can be embedded is called the *page number*. The book-embedding of graphs may be important in several technical applications, e.g., sorting with parallel stacks, fault-tolerant processor arrays design, and layout problems with application to very large scale integration (VLSI). Bernhart and Kainen firstly considered the book-embedding of the planar graph and conjectured that its page number can be made arbitrarily large [JCT, 1979, 320-331]. Heath [FOCS84] found that planar graphs admit a seven-page book embedding. Later, Yannakakis proved that four pages are necessary and sufficient for planar graphs in [STOC86]. Recently, Bekos et al. [STACS14] described an $O(n^2)$ time algorithm of two-page book embedding for 4-planar graphs. In this paper, we embed 5-planar graphs into a book of three pages by an $O(n^2)$ time algorithm.

Keywords: Book embedding; Planar graph; Page number

*Corresponding author. E-mail: ywh222@163.com, yangweihua@tyut.edu.cn

# 1 Introduction

The concept of a book-embedding of a graph was introduced by Ollmann and Kainen. A *book* consists of a line called *spine* and some half-planes called *pages*, sharing the spine as a common boundary. *Book-embedding* of a graph $G = (V, E)$ consists of a (linear) layout $L$ of its notes along the spine $\ell$ of a book(i.e., a one-to-one function from $V$ to $\{1, \ldots, n\}$ ) and the assignment of each edge on the pages so that two edges embedded on the same page do not intersect. We say that two edges $(a, b)$ and $(c, d)$ on the same page cross in the layout $L$ if $L(a) < L(c) < L(b) < L(d)$ or $L(c) < L(a) < L(d) < L(b)$. A central goal in the study of book embedding is to find the minimum number of pages in which a graph can be embedded, this is called the *page number* or *book thickness* of the graph. We denote by $pn(G)$ the page number of $G$. However, determining the page number is a very hard problem. It remains a difficult problem even when the layout is fixed, since determining if a given layout admits a $k$-page book embedding is NP-complete [10].

The book-embedding plays an important role in VLSI design, matrix computation, parallel processing, and permutation sorting. Kapoor et al. in [14] surveyed several applications, as follows (we do not describe the details here): Direct Interconnection Networks [14], Fault-Tolerant Processor Arrays [6], Sorting with Parallel Stacks [7], Single-Row Routing [20], Ordered Sets [17].

The book embedding of graphs has been discussed for a variety of graph families, for examples on complete graphs [3,4], complete bipartite graphs [8,15] and so on (see for examples [13,16,23]). The most famous one of them is the problem determining the page number of planar graphs which has been studied for over 40 years. Bernhart and Kainen [3] firstly characterized the graphs with page number one as the outerplanar graphs and the graphs with page number two as the sub-Hamiltonian planar graphs. They also found that triangulated (maximal) planar graphs requiring at least three pages exist. Moveover, they conjectured that planar graphs have unbounded page number, but this was disproved in [5] and [11]. Buss and Shor [5] investigated a nine-page algorithm based on Whitney's theorem. Heath [11] used a method of "peeling" the graph into levels to reduce the number to seven. Istrail [18] found an algorithm that embeds all planar graphs in six pages. Later, Yannakakis [21] showed that the planar graph admits a four-page book embedding, which can be constructed in linear-time. Since Yannakakis [22] also proved that there are planar graphs that can not be embedded in three pages, four pages are necessary and sufficient for general planar graphs.

One natural question is to consider the page number for specific planar graphs. A planar graph of maximum degree $k$ is called a $k$-planar graph. Heath [12] showed that 3-planar graphs are sub-Hamiltonian, i.e., the 3-planar graph can be embedded into a book of two pages. However, it is NP-complete to decide if a planar graph is sub-Hamiltonian. Bekos et al. [1, 2] described an $O(n^2)$ algorithm to embed the 4-planar graph into a book of two pages, that is, the page number of 4-planar graphs is 2. Bekos et al. also raised a question whether the result can be extended to 5-planar graphs. The page number of 5-planar graphs is clearly at most 4. In this paper, we extend the algorithm in [1, 2] to embed 5-planar graphs into three pages, and thus the page number of 5-planar graphs is either 3 or 2.

The rest of the paper is organized as follows. Section 2 gives an algorithm that every 5-planar graph admits a three-page book embedding and its correctness is proved in Section 3.

## 2 The Algorithm

In this section, we describe an algorithm to embed 5-planar graphs in a three-page book. Similarly to that of Bekos [1, 2], the algorithm is given by a recursive combinatorial construction. It is well known that the page number of a graph equals the maximum of the page number of its biconnected components [3], we therefore assume that $G$ is biconnected. The general idea of our algorithm is as follows. First remove from $G$ cycle $C_{out}$ delimiting the outer boundary of $G$ and contract each bridgeless-subgraph of the remaining graph into a single vertex(say block-vertex). We denote the implied graph by $F$. Note that $F$ is a forest (since $F$ is not necessarily connected). Then, cycle $C_{out}$ is embedded, such that (i) the order of the vertices of $C_{out}$ along the spine $\ell$ is fixed, (ii) all edges of $C_{out}$ are on the same page, except for the edge that connects its outermost vertices. Next, we describe how to assign without crossings: (i) the chords of $C_{out}$, (ii) the edges between $C_{out}$ and $F$, (iii) the forest $F$. Finally, we replace each block-vertex of $F$ with a cycle $C$ delimiting the outer boundary of the bridgeless-subgraph it corresponds to in $G - C_{out}$, and recursively embed its internal vertices and edges similarly to $C_{out}$.

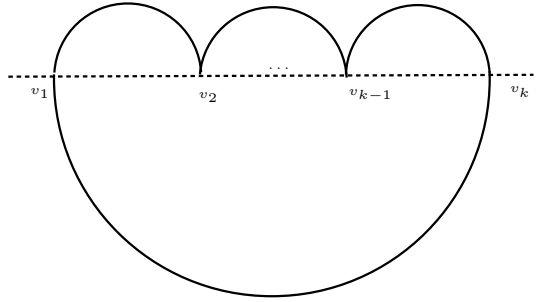

Fig.1

To formalize the idea mentioned above, we consider an arbitrary bridgeless-subgraph, and suppose its outer boundary is a simple cycle $C$. Let $v_1, v_2, \ldots, v_k$ be the vertices of $C$ in the counterclockwise order around $C$. We denote the subgraph of $G$ contained in $C$ by $G_{in}(C)$ and the subgraph of $G$ outside $C$ by $G_{out}(C)$, and let $\overline{G}_{out}(C) = G - G_{in}(C)$ and $\overline{G}_{in}(C) = G - G_{out}(C)$. For the recursive step, we assume the following invariant properties:

IP-1: $\overline{G}_{out}(C)$ has already been completed a three-page book embedding, in which no edge crosses the spine.

IP-2: The combinatorial embedding of $\overline{G}_{out}(C)$ is consistent with a given planar combinatorial embedding of $G$. In other words, the order of the edges in $\overline{G}_{out}(C)$ is consistent with the order of edges in $G$ (the planar embedding) around vertex $v$.

IP-3: The vertices of $C$ are placed in the order $v_1, \ldots, v_k$ along $\ell$, i.e., $L(v_1) < L(v_2) < \ldots < L(v_k)$, and all edges of $C$ are on the same page, except for the edge $(v_1, v_k)$. Without loss of generality, we put edges $(v_i, v_{i+1})$ $(1 \leq i < k)$ of $C$ in the page $p_1$ and edge $(v_1, v_k)$ in page $p_3$ (see Fig.1).

IP-4: If $C$ is not the cycle $C_{out}$ delimiting the outer boundary of $G$, the degree of either $v_1$ or $v_k$ is at most 4 in $\overline{G}_{in}(C)$ since $C$ should be incident to at least one edge of $G_{out}(C)$. Without loss of generality, we assume $deg_{\overline{G}_{in}(C)}(v_k) \leq 4$.

IP-5: If $deg(v_1) = 4$ in $\overline{G}_{in}(C)$, then it is incident to zero or two chords of $C$. If $deg(v_1) = 5$ in $\overline{G}_{in}(C)$, then it is incident to zero or three chords of $C$. This is to assure that $v$ is placed at the right of $v_1$, where $v \in \overline{G}_{in}(C)$.

Note that the combinatorial embedding specified in IP-2 is maintained throughout the whole

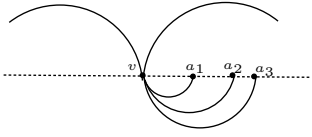
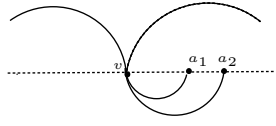
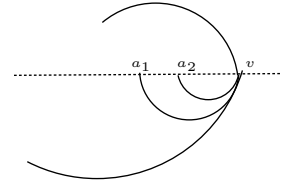

Fig. 2.1  Fig. 2.2  Fig. 2.3

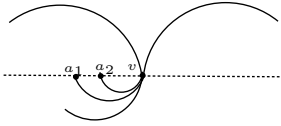
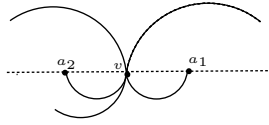
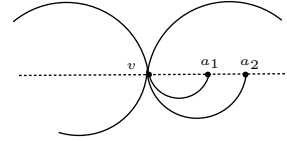

Fig. 2.4  Fig. 2.5  Fig. 2.6

embedding process. Furthermore, it combined with IP-1 is sufficient to ensure planarity.

Again, we describe how to assign the chords of $C$ and the edges between $C$ and $F$. Let $v_i$ be a vertex of $C$, $i = 1, \ldots, k$. Since $G$ is a planar graph with maximum degree five, $v_i$ is incident to at most three non-embedded edges. We denote the edges incident to $v_i$ that follow $(v_i, v_{(i+1) \mod k})$ in the clockwise order of the edges around $v_i$ (as defined by the combinatorial embedding specified by IP-2) by $e_1$, $e_2$ and $e_3$, respectively. Block-vertices that are adjacent to vertices of cycle $C$ are referred to as *anchors*, and block-vertices that are adjacent to other block-vertices only are referred to as *ancillaries*.

Consider an anchor $a$, let $v_{l,a}$ be the leftmost vertex of $C$ adjacent to $a$ along $\ell$. If there is exactly one edge between $a$ and $v_{l,a}$ (i.e., $(a, v_{l,a})$ is simple), we mark this edge. Otherwise, we mark the edge with the largest subscript. Hence, each anchor is incident to exactly one marked edge and each vertex of $C$ is incident to at most three marked edges. Let $v$ be a vertex of $C$, we distinguish four cases for the assignment of the chords and the edges between $C$ and $F$ (i.e., the edges adjacent to $C$).

Case 1. $v$ is adjacent to three anchors $a_1$, $a_2$ and $a_3$ through three marked edges $e_1$, $e_2$ and $e_3$, respectively.

Then, $a_1$, $a_2$ and $a_3$ are placed from left to right and directly to the right of $v$ along $\ell$. Moreover, $e_1$, $e_2$ and $e_3$ are on the page $p_3$ (see Fig.2.1). Note that $v$ cannot be the rightmost vertex of $C$ due to IP-4.

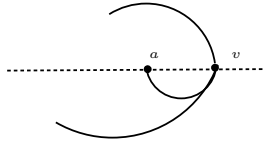
Fig.2.7

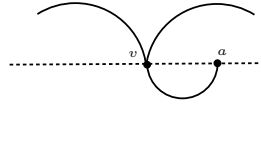
Fig.2.8

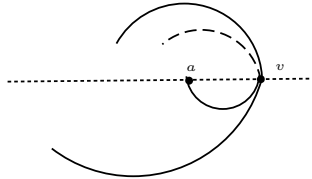
Fig.2.9

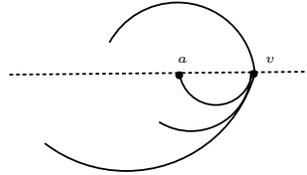
Fig.2.10

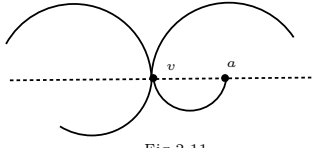
Fig.2.11

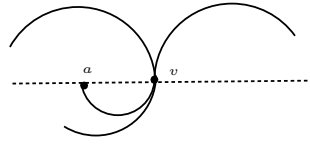
Fig.2.12

Case 2. $v$ is adjacent to two anchors $a_1$ and $a_2$ through two marked edges $e_i$ and $e_j$, respectively.

If $deg(v) = 4$ in $\overline{G}_{in}(C)$, i.e., two marked edges $e_i$ and $e_j$ are edges $e_1$ and $e_2$, respectively. Then $L(a_1) < L(a_2)$ and edges $e_1$ and $e_2$ are on the page $p_3$. Furthermore, we distinguish two sub-cases for the exact placements of $a_1$ and $a_2$. If $v \neq v_k$, then $a_1$ and $a_2$ are placed directly to the right of $v$ (see Fig.2.2); Otherwise $a_1$ and $a_2$ are placed directly to the left of $v$ (see Fig.2.3).

If $deg(v) = 5$ in $\overline{G}_{in}(C)$, $v \neq v_k$ due to degree restriction of $v_k$. In this case, $e_1$, $e_2$ and $e_3$ are placed on the page $p_3$. We distinguish three sub-cases for the exact placements of $a_1$ and $a_2$:

(i) If two marked edges $e_i$ and $e_j$ are edges $e_2$ and $e_3$, $a_1$ and $a_2$ are placed from left to right along $l$ and directly to the left of $v$ (see Fig.2.4);

(ii) If two marked edges $e_i$ and $e_j$ are edges $e_1$ and $e_3$, we place $a_1$ directly to the right of $v$ and $a_2$ directly to the left of $v$ (see Fig.2.5);

(iii) If two marked edges $e_i$ and $e_j$ are edges $e_1$ and $e_2$, $a_1$ and $a_2$ are placed directly to the right of $v$ and $L(a_1) < L(a_2)$ (see Fig.2.6).

Case 3. $v$ is adjacent to one anchor $a$ through the marked edge $e$.

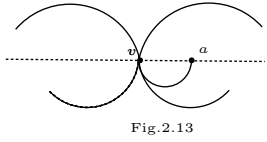
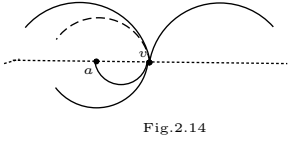
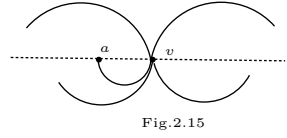

Fig.2.13  Fig.2.14  Fig.2.15

Suppose that $deg(v) = 3$ in $\overline{G}_{in}(C)$, then $e$ is on the page $p_3$. If $v = v_k$, then $a$ is placed directly to the left of $v$ (see Fig.2.7); Otherwise, $a$ is placed directly to the right of $v$ (see Fig.2.8).

Assume $deg(v) = 4$ in $\overline{G}_{in}(C)$. If $v = v_k$, then $a$ is placed directly to the left of $v$. We distinguish two sub-cases for the placements of edges: $e_1$ is marked edge, then $e_1$ is on the page $p_3$ and $e_2$ is on the page $p_2$ (see Fig.2.9); Otherwise, both $e_1$ and $e_2$ are on the page $p_3$ (see Fig.2.10). If $v \neq v_k$, $e_1$ and $e_2$ are placed on the page $p_3$. If $e_1$ is marked edge, $a$ is placed directly to the right of $v$ (see Fig.2.11); Otherwise, $a$ is placed directly to the left of $v$ (see Fig.2.12).

Now we assume $deg(v) = 5$ in $\overline{G}_{in}(C)$. It follows that $v$ is not the rightmost vertex of $C$ by IP-4. In this case, we distinguish three sub-cases.

(i) $e_1$ is the marked edge, then $a$ is placed directly to the right of $v$ and $e_1$, $e_2$ and $e_3$ are on the page $p_3$ (see Fig.2.13).

(ii) $e_2$ is the marked edge, then $a$ is placed directly to the left of $v$. Moreover, $e_3$ is on the page $p_2$ and $e_1$ and $e_2$ are on the page $p_3$ (see Fig.2.14).

(iii) $e_3$ is the marked edge, then $a$ is placed directly to the left of $v$ and all edges incident to $v$ of $C$ are on the page $p_3$ (see Fig.2.15).

Case 4. $v$ is not incident to any marked edge and $deg(v) \neq 2$ in $\overline{G}_{in}(C)$, then the edges incident to $v$ are on the page $p_2$.

Up to now, we have drawn the chords of $C$ and the edges between $C$ and $F$ in pages. We next describe how to embed the edges of $F$. Let $T$ be a tree of the forest of ancillaries (ancillaries form a new forest). We denoted by $\overline{T}$ the tree (*anchored tree*) formed by $T$ and anchors adjacent to some ancillary of $T$.

Suppose that $\overline{T}$ is rooted at the anchor, we assign the vertices of $\overline{T}$ in the order implied by the specific Death First Search (DFS) traversal of $\overline{T}$. If the children $a'$ and $a''$ of $a$ are in the counterclockwise order of edges around $a$, when strating from edge $(a, p(a))$ (where $p(a)$ is the parent of $a$), then $L(a') < L(a'')$.

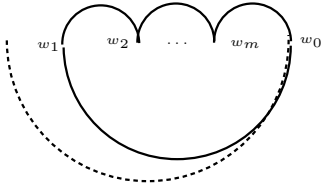 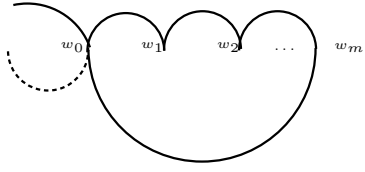 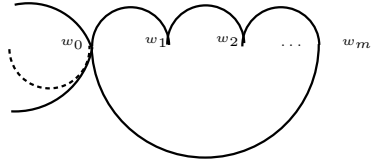

Fig.3.1        Fig.3.2        Fig.3.3

To define the order of the trees in the forest of ancillaries, we create an auxiliary digraph $G_{aux}^T$ whose vertices correspond to trees and there is a directed edge $(v_{\overline{T}_1}, v_{\overline{T}_2})$ in $G_{aux}^T$ if and only if $\overline{T}_1$ has an anchor that is between two consecutive anchors of $\overline{T}_2$.

**Lemma 2.1** ( [1,2]). *Auxiliary digraph $G_{aux}^T$ is a directed acyclic graph.*

Lemma 2.1 implies an embedding of the trees in the order defined by a topological sorting of $G_{aux}^T$. In other words, if $\overline{T}_1$ has an anchor that is between two consecutive anchors of $\overline{T}_2$ along $\ell$, then the tree $\overline{T}_1$ will be embedded before $\overline{T}_2$.

Up to now, $\overline{G}_{in}(C)$ has been embedded, s.t., every bridgeless-subgraph of $G_{in}(C)$ is contracted into a single vertex along $\ell$ and each edge is assigned to one of page $p_1$, page $p_2$ and page $p_3$. Next, we describe how to recursively proceed. Let $a$ be a block-vertex of $G_{in}(C)$ with outer boundary $F_a$: $w_0 \to w_1 \to \ldots \to w_m \to w_0$. We consider two cases:

Case 1. $a$ is an anchor and $w_0$ is incident to a marked edge $e$.

We place $w_0$ as the rightmost vertex of $F_a$ on $\ell$, $w_1$ as the leftmost vertex of $F_a$ on $\ell$, and $w_i$ on the left of $w_{i+1}$ for $i = 1, \ldots, m-1$ and there are no vertices between them (see Fig.3.1). Assign $(w_0, w_1)$ on the page $p_1$ and remaining edges on the page $p_2$. This placement is infeasible when there is an edge or two edges incident to $w_0$ between $(w_0, w_1)$ and the marked edge $e$ in the counterclockwise order of the edges around $w_0$ when starting from $(w_0, w_1)$. In this case, $w_0$ is to the left of $w_1, \ldots, w_m$, i.e. $w_0$ is the leftmost vertex of $F_a$ on $\ell$. This is not possible if there is also an edge $(w_0, w\prime)$ incident to $w_0$ between $(w_0, w_1)$ and the marked edge $e$ in the clockwise order of the edges around $w_0$ when starting from $(w_0, w_1)$. We shall address the problem if $(w_0, w\prime)$ is placed on the page $p_2$.

**Lemma 2.2** ( [1,2]). *Ancillary $a$ can be repositioned on $\ell$ such that: (i) $a$ is placed between two consecutive anchors of T. (ii) The combinatorial embedding specified by IP-2 is preserved and* the edges $(w_0, w)$, $(w_0, w')$ and $(a, w'')$ are on the page $p_2$ and crossing-free. (iii) $w_0$ is leftmost vertex of $F_a$ and $w_i$ to the left of $w_{i+1}$ for $i = 1, \ldots, m-1$; All edges of $F_a$ are on the page $p_3$, except for $(w_0, w_m)$.

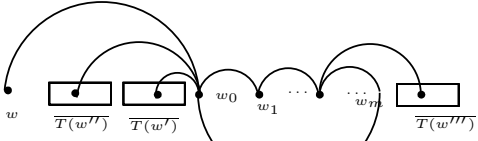
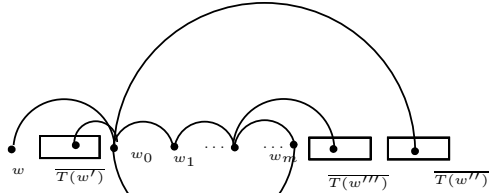

Fig.4.1　　　　　　　　　　Fig.4.2

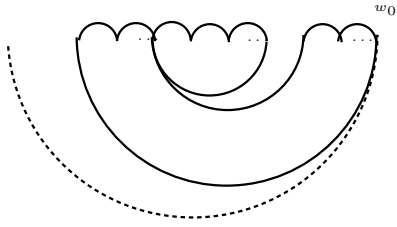
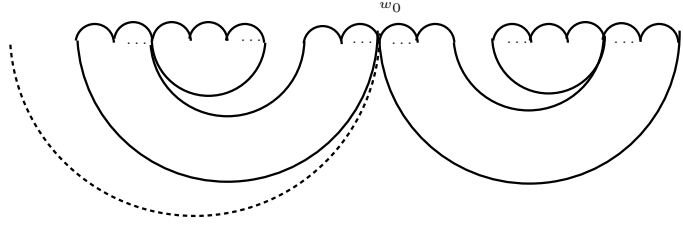

Fig.5.1　　　　　　　　　　Fig.5.2

Case 2. $a$ is an ancillary and $w$ is its parent in the labeled anchored tree $\overline{T}$ in which $a$ belong to.

We place $w_0$ to the leftmost vertex of $F_a$ on $\ell$, $w_1$ to the rightmost vertex of $F_a$ on $\ell$, and $w_i$ to the left of $w_{i+1}$ for $i = 1, \ldots, m-1$ and there are no vertices in between. Moreover, $(w_0, w)$ is on the page $p_1$ and remaining edges are on the page $p_2$. This placement is infeasible in two cases. One case is that there are two edges or only one edge incident to $w_0$ between $(w_0, w_m)$ and $(w_0, w)$ in the counterclockwise order of the edges around $w_0$ when starting from $(w_0, w_m)$, say $(w_0, w')$ (and $(w_0, w'')$). By Lemma 2.2, the vertices of subtrees of $\overline{T}$ rooted at $w'$ and $w''$ are placed to the left of $w_0$ and right of $w$. The other case is that there is an edge $e'$ (and edge $e''$) incident to $w_0$ between $(w_0, w_m)$ and $(w_0, w)$ in the counterclockwise(clockwise, respectively) order when starting from $(w_0, w_m)$. In this case, we place the vertices of subtrees of $\overline{T}$ rooted at $w'$ to right of $w$ and place $w''$ to the left of $w_0$ and the right of $w$.

We assume $F_a$ consisting of simple cycles (each of them is called a subcycle). Note that any two subcycles share at most one vertex of $F_a$ and any vertex of $F_a$ is incident to at most two subcycles. We create a tangency graph $G_{tan}$ whose vertices correspond to subcycles and there is an edge between every pair of subcycles that share a vertex, then $G_{tan}$ is a tree. Let $G_{tan}$ be rooted at the cycle containing $w_0$. Since the degree of $w_0$ is at most 5, $w_0$ lies in at most two subcycles. We firstly suppose that $w_0$ lies in a subcycle, then we place the subcycles of $F_a$ in the order implied by the Breadth First Search (BFS) traversal of $G_{tan}$. Assume $w_0$ lies in two subcycles, say $C_1$ and $C_2$ that correspond to vertices $c_1$ and $c_2$ of $G_{tan}$, respectively. We place the subcycles of $F_a$ to the left (right) of $w_0$ in the order by the Breadth First Search (BFS) traversal of the subtree rooted at $c_1$ ($c_2$ respectively).

Remark. For $\overline{G}_{in}(C)$, the edges assigned to the page $p_1$ are the edges on $C$ (except for $(v_1, v_k)$). The edges on the page $p_2$ consist of the edges on $F$, the edges in Case 4 and the partially edges adjacent to $C$ in Fig.2.9 and Fig.2.14. The edges embedded to the page $p_3$ are edges incident to the vertices of $C$ and $(v_1, v_k)$.

## 3  Proof of Correctness

In this section, we shall show that the algorithm is correct, i.e., no two edges intersect on the same page and IP-1 up to IP-5 are satisfied for arbitrary simple cycle $C$.

**Lemma 3.1** ( [1,2]). *For each ancillary $a$ of a labeled anchored tree $\overline{T}$ there is (i) at least an anchor of $\overline{T}$ with label smaller than that of $a$ and (ii) at least another with label greater than that of $a$.*

**Lemma 3.2.** *If $v \in \overline{G}_{in}(C)$, then $v$ lies on the left of $v_k$ and the right of $v_1$, where $C = \langle v_1, \ldots v_k \rangle$.*

*Proof.* (i) If $v$ is a vertex on $C$, then $v$ must be on the right of $v_1$ and the left of $v_k$ due to IP-3.

(ii) Suppose that $v$ is an anchor. According to the placement of anchors, we have $L(v_1) < L(a) < L(v_k)$ for each anchor $a$ in $\overline{G}_{in}(C)$, i.e., all anchors are placed between $v_1$ and $v_k$.

(iii) Assume that $v$ is an ancillary. Lemma 3.1 implies that there is at least an anchor $a_i$ and

an anchor $a_j$, such that $L(a_i) < L(v) < L(a_j)$. Then $v$ is placed between $v_1$ and $v_k$ according to (ii). □

**Lemma 3.3** ( [1,2]). *The placement of the anchored tree $\overline{T}$ is planar.*

**Lemma 3.4.** *There is no conflict between edges assigned to the pages $p_1$, $p_2$, and $p_3$ in $\overline{G}_{in}(C)$.*

*Proof.* We firstly consider the edges embedded on $p_1$ in $\overline{G}_{in}(C)$: all edges on $C$ except for the one that connects to its outermost vertices. Clearly, the edges assigned on $p_1$ have no crossings.

We now show that all edges assigned to $p_2$ can be embedded in the interior of $\overline{G}_{in}(C)$ without crossings. Lemma 3.3 implies that the edges of $\overline{T}$ has no crossings. Note that there is a path $P$ consisting of edges in $p_2$ joining a pair of consecutive anchors (say $u_1$ and $u_{l+1}$) of $\overline{T}$ and our algorithm must place an ancillary $a$ of $\overline{T}$ between them. Since $c$ is nested by an edge of $P$ and all edges of $\overline{T}$ belong to $p_2$, an edge connecting $a$ to an ancillary of $\overline{T}$ placed between another pair of consecutive anchors of $\overline{T}$ would cross $P$.

By an argument similar to that of Lemma 11 in [1], it can be seen that there is a path $P(u_0 \to u_{l+1})$: $u_0 \to u_{j_1} \ldots u_{j_p} \to u_{l+1}$ consisting of vertices of $\{u_0, \ldots, u_{l+1}\}$, whose edges belong to $p_2$ and for each edge of $P(u_0 \to u_{l+1})$ there is no edges in $p_2$ with endpoints in $\{u_0, \ldots, u_{l+1}\}$ that nests it. Next, $P(u_0 \to u_{l+1})$ contains at least one vertex of $C$. According to the placement of edges, there is only one edge incident to $C$ in $\overline{G}_{in}(C)$, a contradiction. This give a planar embedding of $F$. Hence, there is no crossings between edges assigned in the page $p_2$.

Finally, each edge in $p_3$ is incident to some vertex $v$ on $C$. Note that $G$ is a planar graph and the order of edges incident to $v$ is consistent with the counterclockwise order of edges around $v$ on $C$. Therefore, the edges assigned to page $p_3$ do not cross each other due to IP-2 and IP-1. □

According to Lemma 3.2, the edges of $\overline{G}_{in}(C)$ do not cross the edges of $\overline{G}_{out}(C)$ on the same page. There is no conflict between edges in the same page for $\overline{G}_{out}(C)$ and $\overline{G}_{in}(C)$ due to IP-1 and Lemma 3.4, respectively. Therefore, no edges intersect in the same page. We are now ready to describe that IP-1 up to IP-5 hold when a simple cycle $C_s$ is recursively drawn. Firstly, we prove that IP-1 up to IP-5 are satisfied for $C_{out}$, i.e., the first step of recursion holds.

**Lemma 3.5** ( [1,2]). *Any planar graph $G$ admits a planar drawing $\Gamma(G)$ with a chordless outer boundary.*

IP-1, IP-2 and IP-3 are clearly satisfied. Since $\overline{G}_{out}(C_{out}) = C_{out}$ for $C_{out}$. Lemma 3.5 implies IP-5. If there is a vertex $v \in V(C_{out})$ with $deg_{\overline{G}_{in}(C_{out})}(v) \leq 4$, then it is chosen as $v_k$ and all invariant properties of our algorithm are satisfied. However, if such a vertex does not exist, i.e., each vertex $v \in V(C_{out})$ has $deg_{\overline{G}_{in}(C_{out})}(v) = 5$, this violates IP-4. This case will be addressed by the following lemma.

**Lemma 3.6.** *If $C = C_{out}$, then IP-4 does not necessarily hold.*

*Proof.* If $C = C_{out}$ and IP-4 is not satisfied, then there is no vertex $v$ of $C_{out}$ with $deg_{\overline{G}_{in}(C_{out})}(v) \leq 4$. Suppose $v_k$ is adjacent to vertices $b_1$, $b_2$ and $b_3$ through edges $e_1$, $e_2$ and $e_3$, where $b_1$, $b_2$ and $b_3$ belong to the bridgeless-subgraphs corresponding to anchors $a_1$, $a_2$ and $a_3$ (or simply, $b_1$, $b_2$ and $b_3$ belong to anchors $a_1$, $a_2$ and $a_3$), respectively. The vertices of $\overline{G}_{in}(C)$ should be on the left of $v_m$. We therefore consider the following four cases.

Case 1. All the edges incident to $v_k$ are marked edges.

Case 2. $e_1$ and $e_2$ are marked edges and $e_3$ is non-marked edge.

Case 3. $e_1$ and $e_3$ are marked edges and $e_2$ is non-marked edge.

Case 4. $e_1$ is marked edge and $e_3$ and $e_2$ are non-marked edges.

For Case 1 and Case 2, $b_1$, $b_2$ and $b_3$ belong to the distinct anchors. We augment $G$ by introducing three vertices (say $v_{m+1}$, $v_{m+2}$ and $v_{m+3}$) to the right of $v_m$ (see Fig.6.1). Then IP-4 holds for the augmented graph $G_{aug}$. We claim that $b_1$, $b_2$ and $b_3$ are on the left of $v_m$ in the placement of the augmented graph $G_{aug}$. Let $C_{aug}$ be the outer boundary of the augmented graph $G_{aug}$, and $a_1^{aug}$, $a_2^{aug}$ and $a_3^{aug}$ be the block-vertices belonging to $b_1$, $b_2$ and $b_3$ in $G_{aug} - C_{aug}$, respectively. Since $b_1$ is adjacent to $v_{m+3}$, and $v_{m+3}$ is the rightmost vertex in $G_{aug}$, $b_1$ is placed to the left of $v_m$ according to Fig.2.7. Note that $(v_m, v_{m+2})$ is a chord and the arguments in Fig.2.4 (Fig.2.14) in Case 1 (Case 2, respectively), see Fig.6.2 (Fig.6.4, respectively). We have $b_2$ and $b_3$ are placed to the left of $v_m$. Between $v_m$ and $v_{m+3}$ no vertices of $G_{aug}$ exist, except for $v_{m+1}$ and $v_{m+2}$. Hence, all vertices of $G_{aug} - C_{aug}$ are to the left of $v_m$. If we contract $v_m$, $v_{m+1}$, $v_{m+2}$ and $v_{m+3}$ back into $v_m$, we obtain a valid embedding of $G$ (see Fig. 6.3 and Fig.6.5).

For Case 3 and Case 4, $a_1 \neq a_2$ and $a_1 \neq a_3$, but $a_3 = a_2$ is possible. We augment $G$ by introducing a vertex $v_{m+1}$, such that $v_{m+1}$ is adjacent to $a_3$, $a_2$ and $w_m$ (see Fig.6.6). If $a_3 \neq a_2$, then $v_{m+1}$ belongs to another block-vertex (containing only $v_{m+1}$) and is adjacent to $v_m$. Note that both $a_2^{aug}$ and $a_3^{aug}$ are ancillaries. Since $v_{m+1}$ is adjacent to exactly one vertex of $C_{out}$ (i.e., $v_m$), then, $v_{m+1}$ is to the right of $a_1$ and the left of $v_m$ according to Fig.2.3. Furthermore, $(v_m, v_{m+1})$ is placed on the page $p_3$ and the $(a_2^{aug}, v_m)$ and $(a_3^{aug}, v_m)$ on the page $p_2$ (see Fig.6.7). There is no anchors between $v_{m+1}$ and $v_m$ due to $deg_{G_{aug}}(v_k) = 4$. Therefore, the rightmost anchor of $G_{aug} - C_{aug}$ is $v_{m+1}$. Then, all vertices of $G_{aug} - C_{aug}$ are to the left of $v_{m+1}$. We obtain a valid embedding of $G$ by contracting $v_m$ and $v_{m+1}$ back into $v_m$ (see Fig.6.8).

Suppose $a_3 = a_2 = a$, then $v_{m+1}$ must belong to $a$, and $b_2$, $v_{m+1}$ and $b_3$ appear in the counterclockwise traversal of the outer boundary $C_a$ of $a$. For Case 3, since $a$ is adjacent to $v_m$ through a marked edge $e_3$ in $G$ and through an edge $(v_m, v_{m+1})$ in $G_{aug}$, respectively. Then, $(v_m, v_{m+1})$ is a marked edge, $a$ is directly to the right of $a_1$ and the left of $v_m$, with $v_{m+1}$ being the rightmost vertex of $C_a$ due to Fig.2.3(see Fig.6.9). There is no vertices between $v_{m+1}$ and $v_m$ of $G_{aug}$, since $deg(v_m) = 4$ in $G_{aug}$. So, the rightmost anchor of $G_{aug} - C_{aug}$ has $v_{m+1}$ as its rightmost vertex. Then all vertices of $G_{aug} - C_{aug}$ are to the left of $v_{m+1}$. If we contract vertices $v_m$ and $v_{m+1}$ back to $v_m$, and place $(b_3, v_m)$ on the page $p_2$ and $(b_2, v_m)$ and $(b_1, v_m)$ on the page $p_3$, we obtain a valid drawing of $G$ (see Fig.6.10). For Case 4, $(v_m, v_{m+1})$ is a non-marked edge in $G_{aug}$, since the edges incident to $a$ ($e_2$ and $e_3$) are non-marked edges in $G$. $a$ is directly to the right of $v_i$ ($v_i$ is incident to $a$ thought a marked edge) and the left of $v_m$ due to Fig.2.9 (see Fig.6.11). Hence, we obtain a valid embedding of $G$ by removing $v_{m+1}$ and embedding $(b_3, v_m)$ and $(b_2, v_m)$ on the page $p_2$ (see Fig.6.12). □

**Lemma 3.7** ( [1,2]). *Assume that all trees that precede $T$ in a topological sorting of $G_{aux}^T$ have been drawn on the page $p_2$ without crossings by preserving the combinatorial embedding specified by IP-2. When $\overline{T}$ is drown, the combinatorial embedding specified by IP-2 is also preserved.*

**Lemma 3.8.** *Let $v$ be a vertex of $C$ with degree 3 or 2 in $G_{in}(C)$ that is not the left/right-most vertex of $C$. Let also $v_r$ ($v_l$) be its next neighbor on $C$ to its right (left). Our algorithm places $(v, v_r)$ on the page $p_1$. In fact, it can also be placed on the page $p_3$ without crossings, while the*

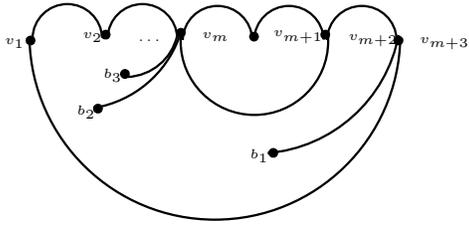

Fig.6.1

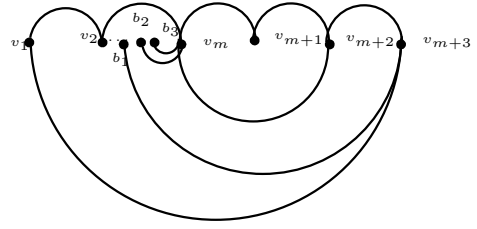

Fig.6.2

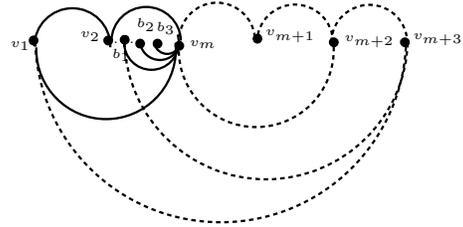

Fig.6.3

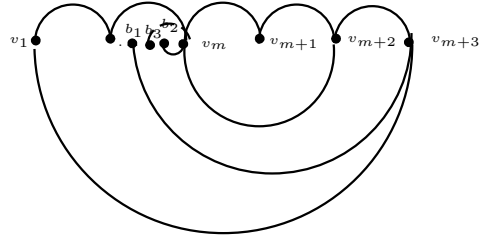

Fig.6.4

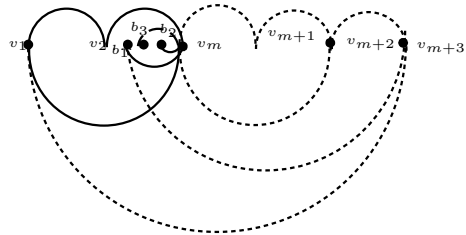

Fig.6.5

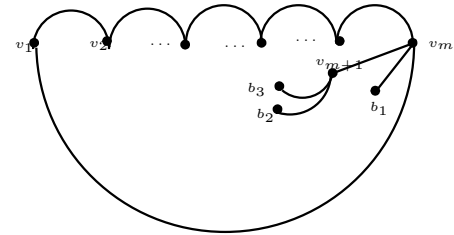

Fig.6.6

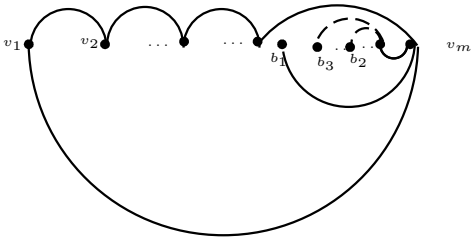

Fig.6.7

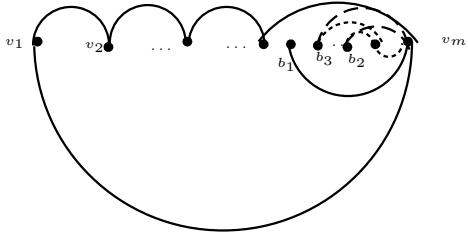

Fig.6.8

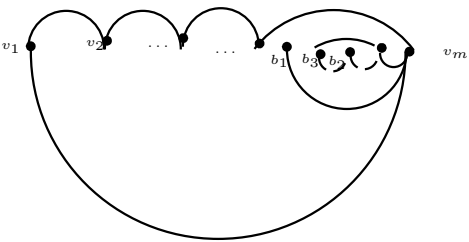

Fig.6.9

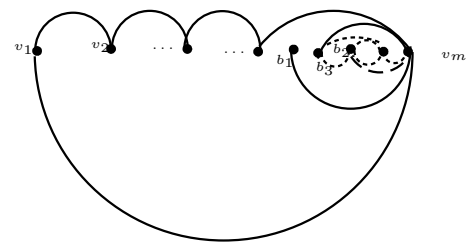

Fig.6.10

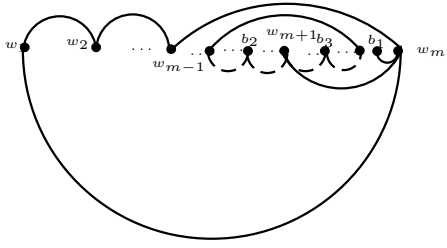
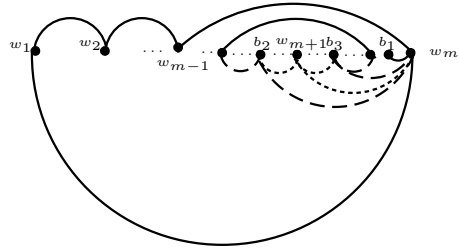

Fig.6.11　　　　　　　　　　　　　　Fig.6.12

*combinatorial embedding specified by IP-2 is maintained.*

*Proof.* Suppose $deg(v) = 2$ in $G_{in}(C)$, then Lemma holds according to Lemma 13 in [1].

Suppose $deg(v) = 3$ in $G_{in}(C)$, if $v$ is incident to non-marked edge $e$, then $(v, v_r)$ can be placed on the page $p_3$ without crossings by the arguments similar to Lemma 13 in [1]. If $v$ is adjacent to an anchor $a$ through a marked edge $e$, then we will move $a$ to the right of $v_r$ and $(v, v_r)$ can be placed on the page $p_3$ without crossings by the proof similar to Lemma 13 in [1]. □

We are ready to describe how the recursive step holds as following: IP-1 up to IP-5 hold for an arbitrary simple cycle $C_s$. Each edge is embedded on one of the three pages: page $p_1$, page $p_2$, and page $p_3$, and no two edges intersect on the same page. Therefore, IP-1 is satisfied. Lemma 3.7 implies IP-2. If $C_s$ is the outer boundary of a block-vertex or a leaf of the tangency tree, then IP-3 trivially holds. If $C_s$ is the outer boundary of a non-leaf of the tangency tree, it contains at least one edge on the page $p_3$. This violates IP-3. In this case, we re-embedding it on the page $p_1$ using Lemma 3.8. IP-4 holds, since suppose that $C_s$ is the outer boundary of a block-vertex or root of the tangency tree of a non-simple outer boundary, then at least one vertex of $C_s$ is adjacent to $G_{out}(C_s)$. Assume that $C_s$ is the outer boundary of internal node of the tangency tree of a non-simple outer boundary, then its leftmost vertex has two edges in $G_{out}(C_s)$. The following lemma implies that IP-5 does not necessarily hold for simple cycle $C_s$.

**Lemma 3.9.** *IP-5 does not necessarily hold for arbitrary cycle $C$.*

*Proof.* If IP-5 does not hold for some simple cycle $C_s$: $w_1 \to w_2 \to \ldots \to w_m$, then we may distinguish three cases.

Case 1. $deg_{\overline{G}_{in}(C)}(w_1) = 4$ and $w_1$ is incident to a chord of $C$.

Case 2. $deg(w_1) = 5$ in $\overline{G}_{in}(C)$ and $w_1$ is incident to a chord of $C$.

Case 3. $deg(w_1) = 5$ in $\overline{G}_{in}(C)$ and $w_1$ is incident to exactly two chords of $C$ (say $(w_1, w_{i_1})$ and $(w_1, w_{i_2})$, where $i_1 \in \{3, \ldots, m-1\}$ and $i_1 < i_2$).

For Case 1 and Case 2, $w_1$ is incident to exactly one chord (say $(w_1, w_i)$, $i \in \{3, \ldots, m-1\}$) of $C$. In general, $(w_1, w_i) \in P(w_1 \to w_j)$, where $P(w_1 \to w_j)$ is a path of chords from $w_1$ to $w_j$ on the page $p_3$. Suppose $w_q \in P(w_1 \to w_j)$, $(w_q, w_x)$ and $(w_q, w_y)$ are chords, $x < y$, we chose $(w_q, w_y) \in P(w_1 \to w_j)$. The restriction implies that $P(w_1 \to w_j)$ is uniquely defined(see Fig.7.1). We refer to it as the separating path of chords of $C_s$, since it splits $\overline{G}_{in}(C_s)$ into two subgraphs: (i) $\overline{G}_{in}(C_l)$ with outer boundary $C_l$ consists of the edges $(w_1, w_2)$, $(w_2, w_3)$, $\ldots$, $(w_{j-1}, w_j)$ and the edges of $P(w_1 \to w_j)$ (see Fig.7.2), and (ii) $\overline{G}_{in}(C_r)$ with outer boundary $C_r$ consists of the edges $(w_j, w_{j+1})$, $\ldots$, $(w_{m-1}, w_m)$, $(w_m, w_1)$ and the edges of $P(w_1 \to w_j)$ (see Fig.7.3).

We now describe how the two sub-instances $\overline{G}_{in}(C_l)$ and $\overline{G}_{in}(C_r)$ can be recursively solved. Since $v$ must be the right of $v_1$, where $v$ is a vertex of $\overline{G}_{in}(C)$. According to the assignment of anchors, we consider the following sub-cases.

Sub-case 1.1: $deg_{\overline{G}_{in}(C)}(w_1) = 4$ and $w_1$ is incident to exactly one chord $e_1$ of $C$.

Sub-case 2.1: $deg(w_1) = 5$ in $\overline{G}_{in}(C)$ and $w_1$ is incident to exactly one chord $e_2$ of $C$.

Sub-case 2.2: $deg(w_1) = 5$ in $\overline{G}_{in}(C)$ and $w_1$ is incident to exactly one chord $e_1$ of $C$.

Observe that if $i \neq j$, then $C_l$ is not simple, i.e., $C_l$ consists of smaller simple subcycles, for which IP-5 holds (hence they can be recursively drawn), except for the first one, that is leftmost embedded along $\ell$. Note that $j \neq m$. If $j = m$ and $deg(w_1) = 3$ ($deg(w_1) = 4$) in $\overline{G}_{in}(C_r)$, then $e_3$ is a bridge according to the placement of chords or $e_2$ is not a chord for some vertex $v_s$($v_s$ belong to $P(w_1 \to w_j)$). If $e_3$ is a bridge, a contradiction since $G$ is biconnected. If $e_2$ is not a chord for some vertex $v_s$, a contradiction since $v_s$ belongs to $P(w_1 \to w_j)$. Since (i) $e_2$ is a marked edge, then $e_3$ is on the page $p_2$ for $v_s$ due to Fig.2.14; (ii) $e_2$ is a non-marked edge, then the edges incident to $v_s$ are placed to page $p_2$ due to Case 4 in placement of edges incident to $C$.

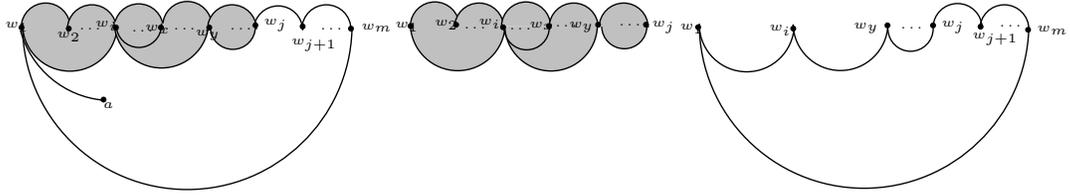

Fig7.1:    Fig7.2: Subgraph $\overline{G}_{in}(C_l)$    Fig7.3: Subgraph $\overline{G}_{in}(C_r)$

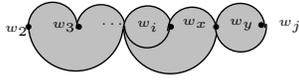 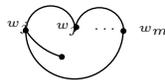 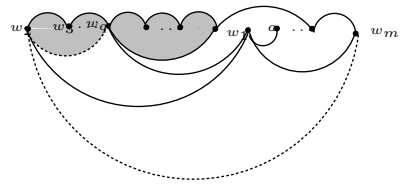

Fig7.4: Subgraph $\overline{G}_{in}(C'_l)$    Fig7.5: Subgraph $\overline{G}_{in}(C'_r)$    Fig7.6:

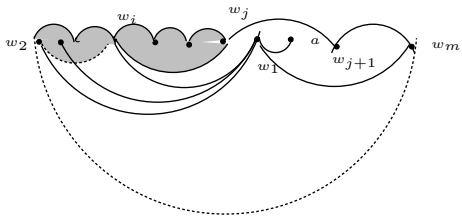 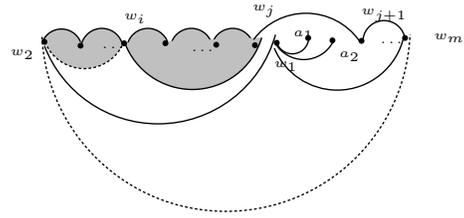

Fig7.7:    Fig7.8:

$\ldots$

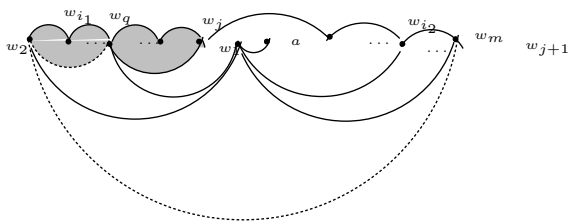 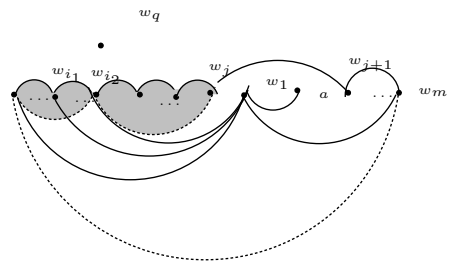

Fig7.9    Fig7.10

In Sub-case 1.1, $deg(w_1) = 2$ in $\overline{G}_{in}(C_l)$, i.e., $deg(w_1) = 3$ in $\overline{G}_{in}(C_r)$. We modify $\overline{G}_{in}(C_l)$ as follows. Remove $w_1$ and join the edges $(w_2, w_i)$ and $(w_2, w_m)$. Then $\overline{G}_{in}(C'_l)$ (see Fig.7.4) has fewer vertices than $\overline{G}_{in}(C)$. We can benefit from this by proceeding recursively, as we initially did with $\overline{G}_{in}(C)$. Eventually, IP-5 should hold for some vertex $w_p$, otherwise a graph with at most 3 vertices on its outer boundary should have a chord, a contradiction. To complete the embedding of $\overline{G}_{in}(C)$, we remove $(w_2, w_i)$ and $(w_2, w_m)$, and connect $w_j$ to its neighbors in $\overline{G}_{in}(C'_r)$ with its copy in $\overline{G}_{in}(C'_l)$ (no crossings are introduced, since the two copies of $w_j$ in $\overline{G}_{in}(C'_l)$ and $\overline{G}_{in}(C'_r)$ are consecutive on $\ell$). It remains to replace the copy of $w_j$ in $\overline{G}_{in}(C'_r)$ with $w_1$, and add $(w_1, w_2)$ and $(w_1, w_i)$ (see Fig.7.6).

In Sub-case 2.1, $deg(w_1) = 3$ in $\overline{G}_{in}(C_l)$ and $deg(w_1) = 3$ in $\overline{G}_{in}(C_r)$. We modify $\overline{G}_{in}(C_l)$ similarly to that of Subcase 1.1, we obtain a valid placement by placing exactly $w_1$ and $a_1$ to the right of $w_j$ and $a_2$ to the right of $w_p$, where $a_1$ and $a_2$ are incident to $w_1$ through marked edges $e_1$ and $e_2$, respectively (see Fig.7.7).

In Sub-case 2.2, $deg(w_1) = 2$ in $\overline{G}_{in}(C_l)$ and $deg(w_1) = 4$ in $\overline{G}_{in}(C_r)$. Similarly to Subcase1.1, $w_1$, $a_1$ and $a_2$ are placed to the right of $w_j$ and the left of $w_{j+1}$, where $a_1$ and $a_2$ are incident to $w_1$ through marked edges $e_1$ and $e_2$, respectively. Hence, we obtain a valid embedding (see Fig.7.8).

For Case 3, we consider two subcases due to the placement of anchors.

Sub-case 3.1: $deg(w_1) = 5$ in $\overline{G}_{in}(C)$ and $w_1$ is incident to exactly two chords $e_1$ and $e_3$ of $C$.

In this case, $P(w_1 \to w_j)$ is defined on the page $p_3$: $(w_1, w_{i_1})$ belongs to $P(w_1 \to w_j)$, if $w_q \in P(w_1 \to w_j)$, $(w_q, w_x)$ and $(w_q, w_y)$ are chords, $x < y$, then $(w_q, w_y) \in P(w_1 \to w_j)$. An analogous argument in Subcase 1.1 implies that $w_1$ and $a$ are placed directly to the right of $w_j$ (see Fig.7.9).

Sub-case 3.2: $deg(w_1) = 5$ in $\overline{G}_{in}(C)$ and $w_1$ is incident to exactly two chords $e_1$ and $e_2$ of $C$.

In this case, we choose $(w_1, w_{i_2}) \in P(w_1 \to w_j)$. If $w_q \in P(w_1 \to w_j)$, $(w_q, w_x)$ and $(w_q, w_y)$ are chords, $x < y$, then $(w_q, w_y) \in P(w_1 \to w_j)$. Hence, $P(w_1 \to w_j)$ is uniquely defined, and $w_1$ and $a$ are placed directly to the right of $w_j$ with the same way of Subcase 1.1 (see Fig.7.10). □

**Theorem 3.10.** *There is a quadratic-time algorithm to construct book embedding for planar graphs of maximum degree 5 on 3 pages.*

*Proof.* Given the planar graph of maximum degree 5 of $n$ vertices. The computation of the bridgeless-subgraphs, the topological sorting of $G_{aux}^T$, BFS-traversals on the tangency trees $G_{tan}$ and DFS-traversals on the anchored tree $\overline{T}$ can be obtained in linear time. Hence the algorithm runs in $O(n^2)$ time. □

# 4 Conclusions and Open Problems

There are planar graphs of maximum degree 8 requiring at least three pages in [3]. This paper has presented an algorithm for embedding a 5-planar graph in three pages. The algorithm can be shown to have time performance $O(n^2)$, so that it is efficient. The natural open problem is whether 2 pages suffices for 5-planar graphs?

# 5 Acknowledgements

The research is supported by NSFC (No.11671296), SRF for ROCS, SEM and Fund Program for the Scientific Activities of Selected Returned Overseas Professionals in Shanxi Province.

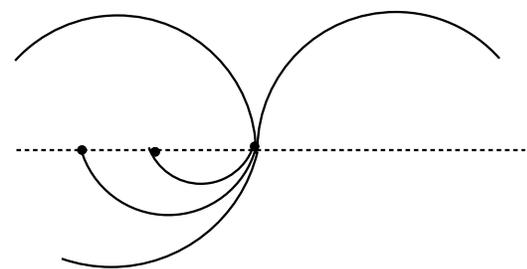 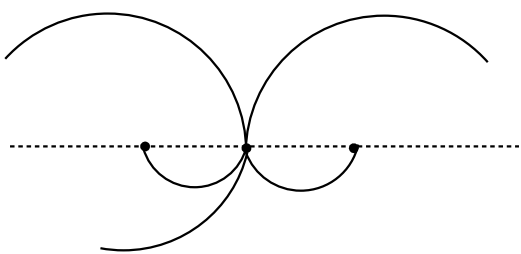 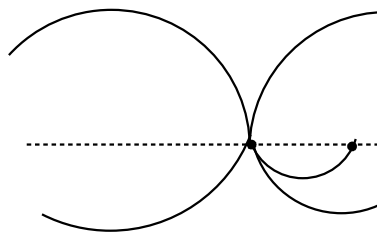